\documentclass[12pt]{article}
\usepackage{amsmath, amssymb, amsthm, amscd}

\numberwithin{equation}{section}

\theoremstyle{plain}
\newtheorem{theorem}{Theorem}
\newtheorem{proposition}{Proposition}
\newtheorem{lemma}{Lemma}

\newtheorem{remark}{Remark}
\theoremstyle{definition}
\newtheorem{definition}{Definition}

\def\ettc{\eta^1_\C \wedge \cdots \wedge \eta_\C^n}

\def\ett{\Omega_s}
\def\bett{\overline \Omega_s}
\def\sett{\sqrt{ \Omega_s}}
\def\ettone{\Omega_1}
\def\bettone{\overline \Omega_1}

\def\C{{\mathbb C}}
\def\R{{\mathbb R}}
\def\h{\hbar}
\def\H{{\cal H}}
\def\L{{\cal L}}
\def\K{{\mathfrak K}}
\def\Hh{{\cal H}^{\rm{\scriptscriptstyle H}}}
\def\dh{\delta^{\rm{\scriptscriptstyle H}}}
\def\Hpr{{\cal H}^{\rm{\scriptscriptstyle prQ}}}
\def\Vpr{{\cal V}^{\rm{\scriptscriptstyle prQ}}}
\def\dpr{\delta^{\rm{\scriptscriptstyle prQ}}}
\def\Hq{{\cal H}^{\rm{\scriptscriptstyle Q}}}
\def\Vq{{\cal V}^{\rm{\scriptscriptstyle Q}}}
\def\dq{\delta^{\rm{\scriptscriptstyle Q}}}
\def\Ha{\rm{\scriptscriptstyle H}}
\def\Qu{\rm{\scriptscriptstyle Q}}
\def\Pr{\rm{\scriptscriptstyle prQ}}

\def\bi{\begin{itemize}}
\def\ei{\end{itemize}}
\def\la{\label}

\newcommand{\be}{\begin{equation}}
\newcommand{\ee}{\end{equation}}
\newcommand{\ba}{\begin{eqnarray}}
\newcommand{\ea}{\end{eqnarray}}

\title{Geometric Quantization, Complex Structures and the Coherent State Transform}
\author{Carlos Florentino$^\dagger$, Pedro Matias$^\ddagger$,
Jos\'e Mour\~ao$^{\dagger}$ \\ and Jo\~ao P.
Nunes$^\dagger$}

\begin{document}

\maketitle

\begin{abstract}
It is shown that the heat operator in the Hall coherent state
transform for a compact Lie group $K$ \cite{Ha1} is related with a Hermitian
connection associated to a natural one-parameter family of complex
structures on $T^*K$. The unitary parallel transport of
this connection establishes the equivalence of (geometric)
quantizations of $T^*K$ for different choices of complex
structures within the given family. In particular,
these results establish a link
between coherent state transforms for Lie groups and
results of Hitchin \cite{Hi} and  Axelrod, Della Pietra and Witten \cite{AdPW}.

\end{abstract}

\newpage

\tableofcontents

\newpage


\section{Introduction}\label{s1}


In the present paper we relate the appearance of the heat operator
in the Hall coherent state transform (CST) for a compact connected
Lie group $K$ \cite{Ha1} with a one-parameter family of complex
structures on the cotangent bundle $T^*K$, in the framework of
geometric quantization.
 The heat equation appears also in the
quantization of $\R^{2n}$ and of Chern-Simons theories
\cite{AdPW,Hi} and in the related theory of theta functions, where
it is associated with the so-called 
Knizhnik-Zamolodchikov-Bernard-Hitchin  (KZBH) connection
\cite{Fa,Las,Ra,FMN}. The general case was studied from a cohomological 
point of view in \cite{Hi}.

Our main motivation is to give a differential geometric interpretation to the appearence
of the heat equation in the K\"ahler quantization of $T^*K$, thus
answering a question raised in \cite{Ha3}, \S 1.3.
This interpretation is based in the projection of the prequantization connection 
to the quantum sub-bundle in geometric quantization, as proposed in \cite{AdPW}. 
The main advantage of this approach consists in the fact that it ensures  
that the quantum connection is Hermitian.
Our method is also complementary
to the one of Thiemann  \cite{Th1,Th2} where
he considers generalized canonical transformations generated
by complex valued functions on the phase space. The heat equation
appears then naturally as the Schr\"odinger equation for these
complex Hamiltonians or complexifiers.  

As shown in \cite{AdPW} \S 1a, for $\R^{2n}$ the
heat equation is associated with independence of the quantization
with respect to the choice of a complex structure within the
family of complex structures which are invariant under
translations.

We consider on $T^*K$  a one-parameter
family of complex structures $\{ J_s \}_{s\in\R_+}$ induced
by the diffeomorphisms
\begin{equation}\label{e11}
\begin{array}{ccccrcl}
T^*K & \simeq & K \times \K^* & \simeq & K\times\K & \stackrel{\psi_s}
{\rightarrow} & K_\C\\
 & & & & (x,Y) & \mapsto & x e^{isY},
\end{array}
\end{equation}
where $K_\C$ is the complexification of $K$.
Here we identify $T^*K$ with $K\times\K^*$ by means of left-translation
and then with $K\times\K$ by means of an $Ad$-invariant inner product
$(\cdot\, ,\cdot )$ on $\K={\rm Lie}(K)$.

Together with the canonical
symplectic structure $\omega$, 
the pair $(\omega,J_s)$ defines on $T^*K$ a K\"ahler structure
for every $s\in\R_+$. Hall has shown \cite{Ha3} that, when one considers
geometric quantization of $T^*K$, the CST, which has been proved to be 
unitary in \cite{Ha1}, gives (up to a constant factor) the pairing 
map between
the vertically polarized Hilbert space and the K\"ahler polarized Hilbert
space, provided that one takes into account the
half-form correction. 

The family of complex structures $\{J_s\}$ is generated by the flow of the 
vector field $v = \sum_j y^j\frac{\partial}{\partial y^j}$. This vector field 
is not Hamiltonian but it is given by $v = J_s (\sum_j sy^j X_j)$, where we note that 
$i \sum_j y^j X_j$
is the Hamiltonian vector field for 
the complex Hamiltonian $\frac{i}{2} |Y|^2$. This is the complex Hamiltonian  
used by Thiemann \cite{Th1,Th2} 
(see also \cite{Ha3}) to generate quantum sates in the holomorphic polarization from 
the vertically polarized ones. The actions of the vector field $v$ and of the Hamiltonian vector field 
corresponding to $s\frac{i}{2}|Y|^2$ coincide on  
$J_s$-holomorphic functions. This explains the relation between 
our formalism and the formalism of complexifiers proposed by Thiemann. 

In order to associate the heat operator with an Hermitian connection,
we collect the prequantum
and quantum Hilbert spaces for all $s\in\R_+$ in a prequantum and
a quantum Hilbert bundles over $\R_+$,
\begin{equation*}
\Hpr \rightarrow \R_+
\end{equation*}
and
\begin{equation}
\nonumber
\Hq \rightarrow \R_+
\end{equation}
and show that the natural Hermitian connection on $\Hpr$ induces
on $\Hq$ a connection given by a heat operator (Theorem \ref{p42}). This connection
turns out to be naturally
equivalent to the connection obtained by varying $\h$ in the
CST of Hall (Theorem \ref{t41}).

Contrary to the flat $\R^{2n}$ case,
and its infinite dimensional generalization considered in \cite{AdPW},
our family of complex
structures is not generated by acting on a fixed one with a family
of canonical transformations. It is generated by the flux of a
vector field which is not symplectic, but rescales $\omega$.
The symplectic structure on $T^*K$ however will be kept fixed
throughout the paper.

Notice that, as could be expected from \cite{Ha3},
the  use of the half-form correction
to define the Hermitian structure on $\Hpr$ plays a decisive
role in the appearance of the heat equation for the
connection induced on the quantum sub-bundle (see also remark \ref{crazyze}).

Our approach to the heat operator and the quantum connection in this 
setting, can also be related to the Blattner-Kostant-Sternberg (BKS) pairing on
the quantum bundle $\Hq$. This will be the theme of a work 
in preparation (\cite{FMMN}).


\section{The quantum connection and the heat equation} \label{s3}
\la{s2.1}

Let $K$ be a compact, connected Lie group.
We will consider first the case when $K$ is semisimple
and will comment briefly on the 
case of compact tori, $K=U(1)^n$, at the end of  section \ref{s4.1}.

We start by recalling from \cite{Ha3,Wo} aspects of the
geometric prequantization of $T^*K$
but with a natural one-parameter family of complex
structures generalizing the fixed complex structure
considered by Hall.

The
prequantum Hilbert bundle $\Hpr$ over this family
is endowed with a natural Hermitian connection, $\dpr$.
The quantum connection $\dq$ induced from $\dpr$
by orthogonal projection on the
quantum Hilbert sub-bundle is then automatically
Hermitian. Our main result in the present section
is Theorem \ref{p42}  in which
we show that $\dq$
corresponds in a precise sense to a family of Laplace operators on $T^*K$.


\subsection{Complex structures
and the prequantum Hilbert bundle}


Consider an
$Ad$-invariant inner product $(\cdot\, ,\cdot )$ on $\K = {\rm Lie}(K)$ and
 $\{X_i\}_{i=1}^{n}$, $n ={\dim K}$, a corresponding
orthonormal basis for $\K$ viewed as the space of left-invariant
vector fields on $K$.
The canonical 1-form on $T^*K$ is given by $\theta=\sum_{i=1}^{n}y^{i}
w^i$ where
$(y^1,\ldots ,y^n)$ are the global coordinates on $\K$
corresponding to the basis $\{X_i\}_{i=1}^n$,
and $\{w^i \}_{i=1}^n$ is the basis of left-invariant
1-forms on $K$ dual to $\{X_i\}_{i=1}^n$, pulled-back to $T^*K$ by the 
canonical projection. The canonical
symplectic 2-form is defined as $\omega=-{\rm d}\theta$. We let
$\epsilon$ denote the Liouville volume form on $T^*K$, given by
\begin{equation}\label{e30}
\epsilon=\frac{1}{n!}\omega^{n} \ .
\end{equation}
Following the geometric
quantization program we consider the trivial complex line bundle $L$
over $T^*K$, $L=T^*K\times \C$, with the trivial Hermitian
structure. Sections of this bundle are thus just functions on
$T^*K$.

Using the diffeomorphisms $\psi_s$ between $T^*K$ and $K_\C$
introduced in (\ref{e11}) we produce a family, parameterized by $s\in\R_+$, of
complex structures $J_s$ on $T^*K$
by pulling back the canonical complex structure $J$
from $K_\C$. Explicitly,
\begin{equation*}
J_s = \psi_{s*}^{-1} \circ J \circ \psi_{s*},
\end{equation*}
where $\psi_{s*}$ denotes the push-forward of the map $\psi_s$.

\begin{proposition}\label{p31}
The pair $(\omega,J_s)$ defines a K\"ahler structure
on $T^*K$ for every $s \in \R_+$, whose K\"ahler potential is
$\kappa_s(x,Y)=s|Y|^2$.
\end{proposition}
\begin{proof}
The family of complex
structures $J_s$  can also be generated by pulling back a fixed
complex structure on $T^*K$ with the family of diffeomorphisms
$\varphi_s (x,Y)=(x,sY)$ since
$J_s=(\psi_1 \circ \varphi_s)^{-1}_* \circ J \circ (\psi_1 \circ \varphi_s)_*
= \varphi_{s*}^{-1} \circ J_1 \circ \varphi_{s*}$.
Therefore
$(\omega, J_s)$ defines a K\"ahler structure on $T^*K$ for
every $s\in\R_+$ if and only if
$((\varphi_s^{-1})^*\omega=\omega/s, J_1)$ defines
a K\"ahler  structure on $T^*K$. This follows from the fact that
$(T^*K,\omega, J_1)$ is a K\"ahler manifold as shown in
\cite{Ha3}. In this reference,
the K\"ahler potential of $(T^*K,\omega, J_1)$
is computed to be $\kappa(x,Y)=|Y|^2$. 
Then, the K\"ahler potential $\kappa_s$ for $(\omega,J_s)$ is 
$\kappa_s=(\varphi_s^*\kappa)/s=s|Y|^2$.
\end{proof}

Let $\tilde{X_j}$, $j=1,...,n$, be the vector fields on $T^*K$
generating the right action of $K$ 
lifted to $T^*K$ and given by
\begin{equation}
\label{lift}
\tilde{X_j}(x,Y)=(X_j, [Y,X_j]).
\end{equation}
Therefore,
\begin{equation*}
\psi_{s*}\tilde{X_j}(x,Y)= X_{j,\C}(xe^{isY}),
\end{equation*}
where $X_{j,\C}$ denotes the natural extension of
$X_j$ from a left-invariant vector field
on $K\subset K_\C$ to  the corresponding  left-invariant vector field on $K_\C$.
Let $\{\tilde w^j\}_{j=1}^{n}$ be the
1-forms defined by $\tilde w^j (\tilde X_k) = \delta^j_k$ and 
$\tilde w^j (J_s \tilde X_k)=0$, for $j,k = 1,...,n$.  
For every $s \in \R_+$, consider also 
the frame of  $J_s$-holomorphic  1-forms given by 
\begin{equation*}
\left\{\tilde \eta^j_{s} =  \tilde w^j - i J_s \tilde w^j \right\}_{j=1}^n \ ,
\end{equation*}
where $(J_s  w)(X)=w(J_s X)$, for a vector field $X$ and a 1-form $w$
on $T^*K$.
Consider also the $J_s$-canonical bundle on $T^*K$ whose sections
are $J_s$-holomorphic $n$-forms with natural Hermitian
structure defined as follows.
For a $J_s$-holomorphic
$n$-form $\alpha_s$, let $|\alpha_s|$ be the unique non-negative 
$C^\infty$ function on $T^*K$ such that $\overline \alpha_s \wedge
\alpha_s = |\alpha_s|^2 b\epsilon$, where $b=(2i)^{n}(-1)^{n(n-1)/2}$.
Following \cite{Ha3} we write
\begin{displaymath}
|\alpha_s|^2 = \frac{\overline \alpha_s \wedge \alpha_s}{b\epsilon}.
\end{displaymath}

A global nowhere vanishing (trivializing) $J_s$-holomorphic
section of the $J_s$-canonical bundle is given by
\begin{equation*}
\Omega_s \equiv \tilde \eta_s^1 \wedge \cdots \wedge \tilde \eta_s^n .
\end{equation*}
Let us introduce the bundle
of half-forms. Let $\delta_s$ denote a square root of the $J_s$-canonical
bundle with a fixed trivializing section whose square
is $\Omega_s$. As
in \cite{Ha3} we denote this section by $\sqrt{\Omega_s}$. Smooth
sections of $L \otimes \delta_s$ are of the form
\begin{equation*}
\sigma_s = f \sqrt{\Omega_s},
\quad f\in C^{\infty}(T^*K).
\end{equation*}
The Hermitian structure on the line bundle $L\otimes \delta_s$ is given by
\begin{equation}\label{e35}
\langle \sigma_s, \tilde \sigma_s \rangle  = \bar f \tilde f
\left(\frac{\bett \wedge \ett} {b\epsilon} \right)^{\frac{1}{2}} = 
\bar f \tilde f|\Omega_s|.
\end{equation}

\begin{definition}\label{d31}
The \emph{prequantum bundle} $\Hpr \rightarrow
\R_+$ is the Hilbert vector bundle with fiber over $s\in \R_+$ given by
\begin{equation*}
\Hpr_s = \overline{\Vpr_s}
\end{equation*}
where the bar denotes norm completion, $\Vpr_s$ is 
\begin{equation*}
\Vpr_s = \left\{\sigma_s \in
\Gamma^{\infty}(L\otimes \delta_s) :\,
||\sigma_s||_s^{\Pr}
\, < \infty \right\},
\end{equation*}
$\Gamma^{\infty}(L\otimes\delta_s)$ denotes the space of $C^\infty$ 
sections of the bundle $L\otimes\delta_s$, and
\begin{equation*}
\langle\sigma_s,\sigma_s\rangle_s^{\Pr} = \int_{T^*K} |\sigma_s|^2 \epsilon.
\end{equation*}
\end{definition}

{}From (\ref{e35}) it is easy to see that 
sections of $\Hpr$ of the
form
\begin{equation}\label{e38}
\frac{f}{\sqrt{|\Omega_s|}}\sett,\quad f \in
L^2(T^*K,\epsilon)
\end{equation}
have $s$-independent norm.

We choose the smooth Hilbert bundle structure on $\Hpr$ as the one 
compatible with the global trivializing map
\ba
\label{trivializ}
\R_+ \times L^2(T^*K,\epsilon) & \to & \Hpr \\
\
(s,f) & \mapsto & \frac{f}{\sqrt{|\Omega_s|}}\sett
\ea

\subsection{The prequantum connection -- $\dpr$}
\la{s2.2}

We now introduce a natural Hermitian connection on $\Hpr$. Before giving its precise
definition we state the following proposition, which is a straightforward consequence
of \cite{Ha2,Ha3}.
Let $\eta (Y)$ be the
$Ad_K$-invariant function defined for $Y$ in a Cartan subalgebra
by the following product over the set $R^+$ of positive roots of $\K$,
\begin{equation}\label{e25}
\eta (Y)=\prod_{\alpha\in R^+} \frac{\sinh \alpha (Y)}{\alpha (Y)}.
\end{equation}
Let $dg$ be the Haar measure on $K_\C$.
We then have,
\begin{proposition} \label{l31}
The following  identities hold:
\begin{enumerate}
\item $|\Omega_s| \equiv \sqrt{\frac{\bett \wedge \ett}{b\epsilon}} =
s^{\frac{n}{2}}\eta(sY)$;
\item $dg_s := (\psi_{s})^{*}(dg) =  s^n \eta^2(sY)\epsilon = |\Omega_s|^2\epsilon$,
\end{enumerate}
where $\eta$ is the function on $T^*K$ defined by equation (\ref{e25}).
\end{proposition}

\begin{proof}
In \cite{Ha2,Ha3} it is shown that
\begin{displaymath}
b\,\eta^2(Y)\epsilon = \bettone \wedge \ettone .
\end{displaymath}
This is exactly the first identity with $s=1$. Recall from Proposition \ref{p31}
that $J_s=\varphi_{s*}^{-1} \circ  J_1 \circ \varphi_{s*}$. This 
implies the equality  $\varphi_s^*(J_1 \beta)
= J_s(
\varphi_s^* \beta)$
for all  1-forms $\beta$ on $T^*K$. Therefore,
\begin{displaymath}
\varphi_s^*(\tilde \eta_1^i) = \varphi_s^*(\tilde w^i) - i \varphi_s^* (J_1
\tilde w^i) = \tilde w^i - i J_s \tilde w^i =\tilde \eta_s^i .
\end{displaymath}
Moreover $\varphi_s^*\epsilon = s^n \epsilon$ and this proves the first
equation. For the second identity, let
\begin{displaymath}
\left\{ \eta^i_\C = w^i_\C - i J w^i_\C \right\}_{i=1}^{n}
\end{displaymath}
be a basis of left $K_\C$-invariant $J$-holomorphic 1-forms on $K_\C$, where $w^i_\C$ is the natural
extension of $w^i$ from $K$ to $K_\C$, obtained by left translations. Then,
the Haar measure on
$K_\C$ is given by
\begin{displaymath}
dg =   \frac{1}{b}\ \overline \Omega_\C \wedge \Omega_\C ,
\end{displaymath}
where $\Omega_\C =  \ettc$.
Since $\psi_s^* \circ J = J_s
\circ \psi_s^*$ for all 1-forms on $K_\C$, and
$\tilde w^i=\psi^*_s(w^i_\C)$, we have
\begin{displaymath}
\psi_s^* (\eta^i_\C) = \psi_s^* (w^i_\C) - i\psi_s^* (J
w^i_\C) = \tilde w^i - i J_s \tilde w^i = \tilde \eta^i_s .
\end{displaymath}
Using the previous result we get the desired identity.
\end{proof}


The trivializing section 
$\sqrt{\Omega_s}$ of the half-form bundle $\delta_s$ is canonical (up to a sign) because
it is obtained from geometric quantization data (the one parameter family of 
K\"ahler polarizations $J_s$ on $T^*K$)
with the only additional structure provided by a fixed
$Ad$-invariant inner product in the Lie algebra $\K$.
This motivates the definition 
of the prequantum connection as the connection induced from the canonical
connection on the trivial bundle by the 
trivialization of $\Hpr$ given in (\ref{trivializ}).


\begin{definition}\label{d22}
The   \emph{prequantum connection} $\dpr$ on $\Hpr$ is the connection
for which sections of the form (\ref{e38}) are horizontal
\begin{equation} \label{e39}
\dpr\left( \frac{f}{\sqrt{|\Omega_s|}} \ \sqrt{\Omega_s}\right) = 0,
\end{equation}
for all $f \in L^2(T^*K,\epsilon)$.
\end{definition}

One can also obtain the same connection using a BKS-type pairing on the
prequantum bundle $\Hpr$. This will be explored in \cite{FMMN}.
Note that the
prequantum connection is Hermitian, that is, it is compatible
with the Hermitian structure on $\Hpr$ in the following sense
\begin{equation}\label{e310}
\frac{d}{d s}\, \langle\sigma,\zeta \rangle^{\Pr} =
\langle\dpr_{\frac{\partial}{\partial s}} \sigma ,\zeta\rangle^{\Pr}
+ \langle\sigma ,\dpr_{\frac{\partial}{\partial s}}\zeta \rangle^{\Pr}
\end{equation}
for all smooth sections 
$\sigma ,\zeta$ of $\Hpr$.


\subsection{The induced quantum connection -- $\dq$}\label{s4}


The K\"ahler polarizations $(\omega,J_s)$ enter already the definition of 
the prequantum Hilbert spaces $\Hpr_s$ through the half-form bundles 
$\delta_s$ and the Hermitian structures (\ref{e35}). To define the fibers 
$\Hq_s$ of the quantum Hilbert sub-bundle $\Hq \subset \Hpr$, one considers polarized, 
or $J_s$-holomorphic, sections of $L \otimes \delta_s$.   

Explicitly, for every $s \in \R_+$, consider the frame of left $K$-invariant vector fields on $T^*K$
\begin{equation}
\nonumber
\left\{Z_{j,s} = \frac{1}{2} (X_j -i J_s X_j)\right\}_{j=1}^n.
\end{equation}
Let the polarizations be given, for every $s\in \R_+$, by
\begin{equation*}
{\cal P}^s_{(x,Y)} = {\rm span}_\C
\left\{\bar Z_{j,s}(x,Y)\right\}_{j=1}^{n} ,
\end{equation*}
where $\bar{Z}_{j,s} = \frac{1}{2}(X_j + i J_s X_j)$.
We use the notation $\bar Z\in {\cal P}^s$ for
$\bar Z_{(x,Y)}\in {\cal P}^s_{(x,Y)}$ for all $(x,Y)\in T^*K$.
Note that these polarizations  $\{{\cal P}^s\}_{s\in \R_+}$
converge, as $s$ tends to zero, to the vertical polarization of $T^*K$, spanned at 
every point by
$\{\frac{\partial}{\partial y^i}\}_{i=1}^n$.

\begin{definition} \label{d41}
The \emph{quantum bundle} $\Hq\rightarrow \R_+$ is the Hilbert
sub-bundle of $\Hpr$ with fiber over $s>0$ given
by
\begin{equation*}
\Hq_s = \overline {\Vq_s}
\end{equation*}
where
\begin{equation*}
\Vq_s = \left\{ f \sqrt{\ett} \in \Vpr_s : \ \nabla_{\bar Z} f = 0
,\ \ \forall \ \bar Z \in {\cal P}^s \right\} ,
\end{equation*}
and $\nabla_{\bar Z} = \bar Z - \frac{1}{i \hbar_0} \theta(\bar Z)$ is
the geometric quantization connection defined on the
trivial bundle
$L$ and
$\hbar_0$ is Planck's constant. We call the solutions of
$\nabla_{\bar Z} f=0$ the \emph{polarized sections of L}.
\end{definition}

\begin{proposition}\label{p41}
For every $s>0$, the ${\cal P}^s$-polarized (or $J_s$-holomorphic) sections of
$L$ are the $C^\infty$ functions $f$ on $T^*K$
of the form
\begin{displaymath}
f=F\, e^{-\frac{s|Y|^2}{2\hbar_0}}
\end{displaymath}
where $F$ is an arbitrary $J_s$-holomorphic function on $T^*K$.
\end{proposition}

\begin{proof}
{}From \cite[Proposition 2.3]{Ha3} the solutions of
$\nabla_{\bar{Z}_{j,s}}f=0$ are
\begin{displaymath}
f=F\, e^{-\kappa_s/2\h_0},
\end{displaymath}
where $F$ is a $J_s$-holomorphic function on $T^*K$, so that the result
follows from proposition \ref{p31}.
\end{proof}

Let us denote by $\langle \cdot ,\cdot \rangle^{\Qu}$ the Hermitian structure
on $\Hq$ inherited from $\Hpr$.

We conclude that
the fibers $\Hq_s$ of $\Hq$ are given by
\begin{equation*}
\Hq_s = \left\{\sigma_s = Fe^{-\frac{s|Y|^2}{2\hbar_0}} \ \sqrt{\Omega_s} , 
\ F \ \hbox{is $J_s$-holomorphic and} \ ||\sigma_s||^Q_s < \infty \right\}  \ .
\end{equation*}

The quantum Hilbert bundle
inherits from $(\Hpr, \dpr)$ a Hermitian connection
$\dq$ which we call the \emph{quantum connection}.
The parallel transport with respect to this connection is
automatically unitary and it establishes the invariance of the
quantization of $T^*K$ with respect to the choice of
polarization within the family $\{{\cal P}^s\}_{s\in \R_+}$.

\begin{definition} \label{d42}
The \emph{quantum connection} $\dq$ is the
Hermitian connection induced on $\Hq$ by the
natural connection $\dpr$ on $\Hpr$
\begin{equation*}
\dq = P \circ \dpr \ ,
\end{equation*}
where $P$ denotes the orthogonal projection  $\Hpr\to\Hq$.
\end{definition}

Below in theorem \ref{p42} we will obtain an explicit expression for $\dq$. 
Consider the second order differential operator
$\Delta_\C^s$
on $T^*K$ given by the pull-back, with respect to $\psi_s$,
of the second order Casimir operator on $K_\C$,

\begin{equation*}
\Delta_\C^s =\sum_{i=1}^n (\tilde{X_i})^2-(J_s \tilde{X_i})^2 .
\end{equation*}
Note that $\Delta_\C^s$ takes $J_s$-holomorphic
functions to $J_s$-holomorphic functions.

Let $\widehat \Delta_\C^s$ be the (unbounded) operator on $\Hq_s$ defined 
on its dense domain by
\be
\la{sun3}
\widehat \Delta_\C^s \left(Fe^{-s|Y|^2/\hbar_0}\sqrt{\Omega_s}\right) =
\Delta_\C^s \left[F \right]e^{-s|Y|^2/\hbar_0}\sqrt{\Omega_s}.
\ee

\begin{theorem}
\label{p42}
Let $F$ be the function on $\R_+ \times T^*K$ obtained, for every $s\in\R_+$, 
as the pull-back of a given $J$-holomorphic function $\hat F$ on $K_\C$, 
\be
\label{ts}
F(s,x,Y)= \hat F (xe^{isY}),
\ee   
such that 
$
F(s,\cdot)e^{-\frac{s|Y|^2}{2\hbar_0}} \sett \in \Hq_s
$ 
is in the domain of the operator $\widehat \Delta^s_\C$.
The quantum connection $\dq$ acts on sections of $\Hq$ of
the form 
\be
\label{ts2}
F\, e^{-\frac{s|Y|^2}{2\hbar_0}} \sett
\ee 
as
\begin{equation}
\dq_{\frac{\partial}{\partial s}}[ F\, e^{-\frac{s|Y|^2}{2\hbar_0}} \sett ]
\label{e44}
= \frac{\hbar_0}2 \left(-\frac{1}{2} \Delta_\C^s + |\rho|^2 \right)
\left[F\right]\, e^{-\frac{s|Y|^2}{2\hbar_0}} \sett ,
\end{equation}
where $\rho$ is half the sum of the positive roots of $\K$.
\end{theorem}

Notice that choosing the sections of $\Hq$ in the form (\ref{ts})
and (\ref{ts2}) corresponds to choosing moving frames, more precisely
a class of global moving frames related by $s$-independent transformations. 
Having made such a choice the 
covariant derivative in the direction of $\partial/\partial s$
is defined by a linear operator acting on the fibers as in (\ref{e44}). 

We will divide the proof of this theorem in several lemmata.
Let us denote by $W$ the following vector field on $T^*K$,
$$
W= i\sum_{j=1}^{n} y^j Z_{j,s} \ .
$$
Note that, from (\ref{lift}),
\begin{equation*}
 W =\frac i{2}\sum_{j=1}^{n} y^j 
( X_{j} - i J_s  X_j)=\frac i{2}\sum_{j=1}^{n} y^j 
(\tilde X_{j} - i J_s \tilde X_j).
\end{equation*}

\begin{lemma}
\label{p43} Let $\hat F$ be a fixed $C^\infty$ function on $K_\C$,
$F$ be the function on $\R_+\times T^*K$, given by
$F(s,\cdot) = \psi_s^* \hat F$, 
and let $f\in C^\infty(T^*K)$ be  a function  only of $Y$. We have,
\begin{itemize}
\item[i)] If $\hat F$ is holomorphic then 
\begin{equation*}
\frac{\partial
F}{\partial s} = W  F
\end{equation*}

\item[ii)] If $\hat F$ is right $K$-invariant then 
\begin{equation*}
\frac{\partial F}{\partial s} = 2 \ W  F = 2 \ 
\bar{W} F
\end{equation*}

\item[iii)] 
\begin{equation*}
Wf = \frac {1}{2s} \sum_{j=1}^n y^j 
\frac {\partial}{\partial y^j} \ f.
\end{equation*}

\end{itemize}
\end{lemma}

\begin{proof}
A direct computation gives, 
\begin{equation*}
\frac{\partial}{\partial s}(\psi_s^* \hat F) (x,Y)=
\frac{\partial}{\partial s} \hat F (x e^{isY}) = \sum_{j=1}^n  
 \ y^j J_s \tilde X_j F. 
\end{equation*}
The special cases {\it i)} and {\it ii)}
above follow from this equation.

The identity in {\it iii)} follows from
\begin{equation*}
\sum_{j=1}^{n} y^j 
  J_s  X_j =  \frac 1s \sum_{j=1}^{n} y^j 
   \frac{\partial}{\partial y^j} .
\end{equation*}
 
\end{proof}

\begin{lemma}\label{byparts}
Let $X$ be a smooth vector field on
$T^*K$, with $|X(y^i)| < c\exp(\alpha |Y|)$ for some positive constants $c$ and $\alpha$,
and let $\phi \in C^\infty (T^*K)$ be such that 
\be
\nonumber
\label{uu1}
|\phi(x,Y)| < e^{-\delta |Y|^2},
\ee
for $|Y|>R$ and fixed positive constants $R,\delta$. 
Then, 
\begin{equation*}
\int_{T^*K} {\cal L}_X(\phi\epsilon) =0.
\end{equation*}
\end{lemma}
\begin{proof}
Using Cartan's formula, ${\cal L}_X=d\circ\iota_X+\iota_X\circ d$,
we have ${\cal L}_X(\phi\epsilon) = d(\phi\, \iota_X \epsilon)$, 
where the symplectic volume form is  $\epsilon = w^1 \wedge  dy^1 \wedge 
\cdots \wedge w^n \wedge dy^n$. Since $T^*K \cong K\times \K$,
using Stokes formula we only need to show that  
$$
\lim_{R\rightarrow \infty} \int_{K\times S^{n-1}_R} \phi \,\iota_X \epsilon =0,
$$
where $S^{n-1}_R$ denotes the sphere of radius $R$ in $\K$ centered at the origin.
We have 
$$
\left|\int_{K\times S^{n-1}_R} \phi \,\iota_X \epsilon\right| < c' 
e^{-\frac{\delta}{2}R^2},
$$ 
for some positive constant $c'$, which proves the lemma. 
\end{proof}

Let $\cal B$ be the subspace of the space ${\cal H}(K_{\C})$ 
of holomorphic functions on $K_\C$,
spanned by the holomorphic functions 
\be
\la{sun2} {\rm tr}\,(\pi(g)A),
\ee 
where $\pi$ is an irreducible 
finite-dimensional representation of $K$ 
extended to an holomorphic representation
of $K_\C$
and $A\in {\rm End}\, V_\pi$. 
Let ${\cal F}_s$ denote the subspace of $\Hq_s$ given by sections of the form 
(\ref{ts}) and (\ref{ts2}) with $\hat F\in \cal B$. 
It follows from the Lemma 10 of \cite{Ha1} that ${\cal F}_s$ is dense in $\Hq_s$.

\begin{lemma}\label{p44}
Let $\sigma,\zeta\in\Gamma(\Hq)$ with $\sigma_s,\zeta_s \in {\cal F}_s$, 
$\forall s\in\R_+$,  and
\begin{eqnarray*}
\sigma  = F\, e^{-\frac{s|Y|^2}{2\hbar_0}} \sett ,& &
\zeta  =  G\, e^{-\frac{s|Y|^2}{2\hbar_0}} \sett.
\end{eqnarray*}
Then we have,
\begin{displaymath}
\int_{T^*K} \left( \bar{W} \,
-\, \frac{|Y|^2}{\h_0}+\frac{1}{2}
\frac{\partial\ln|\Omega_s|}{\partial s} + \frac n{4s}\right) \left[\bar{F}\right] G 
e^{-\frac{s|Y|^2}{\h_0}}|\Omega_s|\,\epsilon = 0.
\end{displaymath}
\end{lemma}

\begin{proof}
The vector field $W$ satisfies $\bar W(y^j)=\frac{1}{2s}y^j$ and  we can apply lemma \ref{byparts}
with $0<\delta <\frac{s}{\hbar_0}$ to the first term above.
Integrating by parts we obtain
\begin{eqnarray}\label{e410}
-\int_{T^*K}  \bar{W} \left[\bar{F}\right] G 
e^{-\frac{s|Y|^2}{\h_0}}|\Omega_s|\,\epsilon \!\!\!
& = & \!\!\! \int_{T^*K}\bar{F}\, G\, \bar{W}\,\left[
e^{-s|Y|^2/\h_0}\right]\,|\Omega_s|\,\epsilon \nonumber\\
   & + &\!\!\!\!\!\!\int_{T^*K}\bar{F}\, G\,
e^{-s|Y|^2/\h_0}\,\bar{W}\,\left[\ln|\Omega_s|\right]\,|\Omega_s|\,\epsilon \nonumber\\
   & + & \!\!\!\!\!\!\int_{T^*K}\bar{F}\, G\, e^{-s|Y|^2/\h_0}\,|\Omega_s|\,
\L_{\bar{W}}\,(\epsilon).
\end{eqnarray}
For the first term on the r.h.s. of (\ref{e410}) we have from {\it iii)} in lemma
\ref{p43}
\be 
\la{cc0cc}
\bar
W[e^{-s|Y|^2/\h_0}]=\frac {1}{2s}\sum_{j=1}^n
y^j\frac{\partial}{\partial y^j}e^{-s|Y|^2/\h_0} = -\frac
{|Y|^2}{\hbar_0} e^{-s|Y|^2/\h_0}. \ee
 {}From {\it ii)} in lemma
\ref{p43}, we see that
\be
\la{cc1cc}
\bar W [\ln|\Omega_s|] = \frac{1}{2} \frac{\partial \ln|\Omega_s|}{\partial s}
-\frac{n}{4s}\ .
\ee

{}For the third term in \eqref{e410}, we have 
\begin{displaymath}
\L_{\bar{W}}\,(\epsilon)=\frac{n}{n!}\,\L_{\bar{W}}\,
(\omega)\wedge\omega^{n-1}.
\end{displaymath}
{}From $\omega=-d\theta$, Cartan's formula $\L_{\bar
W}=d\circ\iota_{\bar W}+\iota_{\bar W}\circ d$ and
\begin{displaymath}
\iota_{\bar W}\,d\theta=\frac{1}{2s}\,\theta + \frac{i}{4}\,d|Y|^2 \ ,
\end{displaymath}
we obtain $\L_{\bar W}\,(\omega)=(1/2s)\omega$, which implies that
\begin{equation}\label{e413}
\L_{\bar W}\,(\epsilon)=\frac{n}{2s}\,\epsilon.
\end{equation}
Substituting \eqref{cc0cc}, \eqref{cc1cc} and \eqref{e413} into
\eqref{e410} we obtain the desired result.
\end{proof}

Recall that the prequantum
connection is Hermitian, so it satisfies equation (\ref{e310}).
Consider sections of $\Hq$, 
\be
\label{que}
\sigma = F\, e^{-\frac{s|Y|^2}{2\hbar_0}} \sett,\,\,\qquad \zeta= G\, e^{-\frac{s|Y|^2}{2\hbar_0}} \sett,
\ee
where $F,G$ are as in (\ref{ts}). 
These sections also satisfy (\ref{e310}) and we have
\begin{equation}\label{e43}
\langle \dq_{\frac{\partial}{\partial s}}\sigma,\zeta\rangle^{\Qu} =
\langle \dpr_{\frac{\partial}{\partial s}}\sigma, \zeta\rangle^{\Pr}.
\end{equation}

\begin{lemma} \label{ll1}
Let $\sigma,\zeta$ be as above.
Then, the following identity holds
\begin{eqnarray}\label{e414}
\langle \dq_{\frac{\partial}{\partial s}}\sigma ,\zeta\rangle^{\Qu} =
 s^{-n/2} \int_{K_\C}\bar{\hat F}\hat
G \ \hbar_0 \left(\frac{|Y|^2}{2\h^2}-\frac{n}{4\hbar}\right)
\frac{e^{-|Y|^2/\h}}{\eta(Y)} \ dg  ,
\end{eqnarray}
where $dg$ is the Haar measure on $K_\C$, $g=xe^{iY}$ and $\hbar=s\hbar_0$.
\end{lemma}

\begin{proof}
From the definition of horizontal sections
(\ref{e39}) and from (\ref{e43}), we obtain 
\begin{equation*}
\langle \dq_{\frac{\partial}{\partial s}}\sigma ,\zeta \rangle^{\Qu} = 
 \int_{T^{*}K}
\overline{\left( \frac{\partial F}{\partial s} - \frac{|Y|^2 F}
{2\h_0}+\frac{F}{2}\frac{\partial \ln |\Omega_s|}{\partial s} \right)} \, 
G\, e^{-\frac{s|Y|^2}
{\h_0}}|\Omega_s|\,\epsilon .
\end{equation*}
We now use lemmata \ref{p43} and \ref{p44} to simplify the equation above for
the quantum connection to get,
\[
\langle \dq_{\frac{\partial}{\partial s}}\sigma ,\zeta \rangle^{\Qu}=
\int_{T^*K}\bar
FG\left(\frac{|Y|^2}{2\h_0}-\frac{n}{4s}\right)
e^{-s|Y|^2/\h_0}|\Omega_s| \ \epsilon ,
\]
which with the help of $\psi_s$ and proposition \ref{l31} 
gives (\ref{e414}).\end{proof}

Let us  introduce the $K$-averaged heat kernel measure $d\nu_\h$ on
$K_\C$ given by \cite{Ha3}
\begin{equation}\label{e24}
d\nu_\h (g)=\nu_\h (g) \, dg= c_\h\,\frac{e^{-|Y|^2/\h}}{\eta (Y)}\, dg,
\end{equation}
 and $c_\h \, = (\pi\h)^{-n/2} e^{-|\rho|^2 \h}$, $\rho$ being
half the sum of the positive roots. Recall from \cite{Ha1} that
$\nu_\h$ satisfies the equation
\be \la{00712}
\frac{\partial \nu_\h}{\partial\h} = -\frac{1}{4}\Delta_\C \nu_\h
\ee
on $K_\C$, where
\begin{equation*}
\Delta_\C = \sum_{i=1}^n (X_{i,\C})^2 - (J X_{i,\C})^2 \ ,
\end{equation*}
is the Casimir operator for $K_\C$. The equation (\ref{00712})
is equivalent to the following equality:
\begin{lemma}
\la{ppp1} \be \nonumber
\left(\frac
{|Y|^2}{2\hbar^2}-\frac{n}{4\hbar}\right)\frac{e^{-|Y|^2/\h}}{\eta(Y)}
=
\left(-\frac
{1}{8}\Delta_\C+\frac{|\rho|^2}2\right)
\frac{e^{-|Y|^2/\h}}{\eta(Y)}.
\ee
\hfill $\Box$
\end{lemma}

\bigskip
We are now ready to prove Theorem \ref{p42}.

\begin{proof}
Consider $\sigma$ and $\zeta$ as in (\ref{que})
and let
\be
\nonumber
Z_j = \frac 12 \left(X_{j,\C} -iJ X_{j,\C}\right), \qquad j=1, \cdots  n, 
\ee
so that 
$$
\Delta_\C = 2 \sum_{j=1}^n Z_j^2 + \bar Z_j^2.
$$
{}From (\ref{e414}) and Lemma \ref{ppp1}, since 
$$X_{j,\C}(\frac{e^{-|Y|^2/\h}}{\eta(Y)}) = 0, \quad \hbox{for all} \ j,$$ 
we obtain 
\ba 
\nonumber
 \langle \delta^Q_{\frac \partial {\partial s}}\sigma ,
\zeta\rangle^Q &=& s^{-n/2}\int_{K_\C} \bar {\hat F}\hat G \ \hbar_0 \left(-\frac
{1}{8}\Delta_\C+\frac{|\rho|^2}2\right)
\left[\frac{e^{-|Y|^2/\h}}{\eta(Y)}\right]dg=\\
\nonumber
&=&s^{-n/2}\int_{K_\C} \bar {\hat F} \hat G \ 
\frac{\hbar_0}2 \left(-\sum_{j=1}^n\bar Z_j^2+|\rho|^2\right)
\left[\frac{e^{-|Y|^2/\h}}{\eta(Y)}\right]dg . \qquad\ea 
Lemma \ref{byparts} can be applied since $\bar Z_j (y^j)$ grows 
linearly with $|Y|$ (see section 6 in  \cite{Ha3}).
Using it to integrate twice by parts
with respect to $\bar Z_j$ and noticing that, from the
bi-invariance of $dg$, ${\cal L}_{\bar Z_j}dg=0$ and also $\bar
Z_j(\hat G)=0$ we obtain the statement of the theorem
\be
\label{ts999}
\langle \delta^Q_{\frac \partial {\partial s}}\sigma ,
\zeta\rangle^Q =
\int_{T^*K} \frac{\hbar_0}2 \,\overline{(- \frac 12\Delta_\C^s + |\rho|^2) [F]}\,  G\,  
 e^{-s|Y|^2/\hbar_0}\,|\Omega_s|\epsilon, 
\ee
for sections $\sigma, \zeta$ with values in ${\cal F}_s \ni \sigma_s, \zeta_s$, for all
$s\in \R_+$. 
The operator $\widehat \Delta_\C^s$ in (\ref{sun3})
 is essentially self-adjoint  
(it has a basis of eigenvectors, with $\hat F$'s given by matrix elements
of finite dimensional holomorphic representations as in (\ref{sun2}), 
with real, and nonpositive, eigenvalues) and 
the space ${\cal F}_s$ is dense in $\Hq_s$. Therefore, the expression (\ref{ts999}) 
implies that  $\delta^Q_{\frac{\partial}{\partial s}}$ is given by (\ref{e44}) 
for  sections of $\Hq$ in the form (\ref{ts}) and (\ref{ts2}) 
and with values in the domain of  $\widehat \Delta_\C^s$,
 for all
$s\in \R_+$.
\end{proof}

\begin{remark}\label{crazyze}
Note that the simple expression (\ref{e414}) would not be valid 
without the inclusion of the half-form correction in the definition of 
the Hermitian structure on $\Hpr$ (see (\ref{e35}) and proposition
\ref{l31}). The half-form leads to the 
cancellation of the term proportional to $\partial \ln |\Omega_s|/\partial s$.
This is what ultimately leads to the heat operator in the 
quantum connection.
\end{remark}


\subsection{The heat equation}
\la{s4.1}

Let $\Psi $ be the diffeomorphism
\begin{equation*}
\begin{array}{rll}
\Psi :\R_+\times T^{*}K & \rightarrow  & \R_+\times K_{\C}
\\
(s,(x,Y)) & \mapsto  & (s,xe^{isY}),
\end{array}
\end{equation*}
so that $\Psi (s,.)=\psi _{s},$ for all $s\in \R_+.$ Consider 
sections $\sigma$ of $\mathcal{H}^{Q}$ of the form
\begin{equation}
\sigma_s =\frac{1}{\sqrt{a_{s}}}F(s,\cdot)\ e^{-\frac{s|Y|^{2}}{2\hslash _{0}}}\sqrt{%
\Omega _{s}},  \label{sal2a}
\end{equation}
where $a_s=(\pi\hslash_0)^{n/2}e^{|\rho|^2\hslash_0 s}$,
\begin{equation}
F=\Psi ^{*}\widetilde{F}  \label{sal2b}
\end{equation}
and $\widetilde{F}$ is a $C^{\infty }$ function on $\R_+\times K_{%
\C},$ holomorphic in the second variable. We then have,

\begin{theorem}
\label{th3.1}
A section of $\mathcal{H}^{Q}$ of the form (\ref{sal2a}), (\ref{sal2b}) 
with $\sigma_s$ in the domain of $\widehat \Delta_\C^s$, for all $s\in \R_+$,
is $\delta ^{Q}$-horizontal if and only if 
$\widetilde{F}$ is a solution of the
following heat equation on $K_{\C},$%
\begin{equation*}
\frac{\partial }{\partial s}\widetilde{F}=\frac{\hslash _{0}}{4}\Delta _{%
\C}\widetilde{F} . 
\end{equation*}
\end{theorem}

\begin{proof}
Substituting (\ref{sal2a}) in (\ref{e44}) we obtain
\[
\delta^Q _{\frac{\partial }{\partial s}}\sigma = \Psi^*\left( \frac{\partial }{
\partial s}\widetilde F-\frac{\hslash _{0}}{4}\Delta _{\C}\widetilde F\right)
\frac{e^{-\frac{s|Y|^{2}}{2\hslash _{0}}}}
{\sqrt{a_{s}}}\sqrt{\Omega _{s}}.
\]
\end{proof}

Besides establishing the relation of the heat equation on complex semisimple
Lie groups $K_{\C}$ with the geometric quantization of $T^{*}K$ this result
also provides the link with the coherent state transform introduced by Hall
in \cite{Ha1}.

Concerning the  case 
 $K=U(1)^n$, note that the results in 
proposition \ref{l31} and equations (\ref{e24}) and (\ref{00712}) are still valid 
 if we replace $\eta (sY)$ by $1$ and the Weyl vector $\rho$ by $0$. 
Therefore, all the results above apply also to this case.


\section{The quantum connection and the CST}
\label{s2}

The coherent state transform for $K$ defines a parallel transport on 
a Hilbert bundle over $\R_+$ for an Hermitian connection that we
denote by $\delta^H$.  
In this section we show that this parallel transport
is naturally equivalent
to the parallel transport of the quantum
connection $\dq$. 
This follows from
propositions \ref{l31} and \ref{p41}, and from (\ref{e35}) and (\ref{e24}).

\subsection{The CST connection - $\dh$}

Let $\rho_{\h}$, $\h >0$, be
the heat kernel for the Laplacian $\Delta$ on $K$ associated to the
$Ad$-invariant inner product $(\cdot\, ,\cdot )$ on $\K$,
 $\Delta = \sum_{i=1}^{n} X_i^2$.
As proved in \cite{Ha1}, $\rho_{\h}$ has a unique analytic
continuation to $K_{\C}$, also denoted by $\rho_{\h}$. The
$K$-averaged coherent state transform (CST) is defined as the map
\begin{eqnarray} \label{e21}
C_{\h} & : & L^{2}(K,dx)\rightarrow {\cal H}(K_{\C}) \nonumber \\
(C_{\h}f)(g) & = & \int_{K} f(x) \rho_{\h} (x^{-1} g)\, dx, \qquad
f\in L^{2}(K,dx) , \ g\in K_{\C}  ,
\end{eqnarray}
where $dx$ is the normalized Haar measure on $K$. For each $f\in L^{2}(K,dx)$,
$C_{\h}f$ is the analytic continuation to $K_{\C}$ of the solution
of the heat equation on $K$,
\begin{equation*}
\frac{\partial  u}{\partial \h} =\frac{1}{2} \Delta u,
\end{equation*}
with initial condition given by $u(0,x)=f(x)$. Therefore, $C_{\h}f$ is given by
\begin{equation*}
(C_{\h}f)(g) = ({\cal C}\circ \rho_{\h} \star f)(g) = \left({\cal
C}\circ e^{\frac{\h\Delta}{2}}f \right)(g),
\end{equation*}
where $\star$ denotes the convolution in $K$ and ${\cal C}$ denotes analytic
continuation from $K$ to $K_\C$.
 Recall that the $K$-averaged heat kernel measure $d\nu_\h$ on
$K_\C$ is given by (\ref{e24}).
Hall proves the following:
\begin{theorem}[Hall] \label{t21}
For each $\h> 0$, the mapping $C_\h$ defined in (\ref{e21}) is an unitary
isomorphism from $L^{2}(K,dx)$ onto the Hilbert space $\H L^2
(K_\C , d\nu_\h) :=L^{2}(K_\C,d\nu_\h)\bigcap \H (K_{\C})$.
\end{theorem}

This transformation defines a Hilbert vector bundle
$\Hh$ over the one-dimensional base $\R_+$ with
global coordinate $\h$,
\begin{equation*}
\begin{array}{rcl}
\Hh & \rightarrow & \R_+ \\
\Hh_{\h} & = & \H L^2(K_\C ,d\nu_\h).
\end{array}
\end{equation*}
which we call the \emph{CST bundle}.
The  Hilbert space structure on the
fibers of $\Hh$ is defined by
\be
\nonumber
\langle  F_\h , G_\h \rangle_\h^{\Ha} :=\int_{K_\C}
\overline{ F_\h(g)}\, G_\h(g)\, d\nu_\h(g), \quad
 F_\h, G_h \in \Hh_\h .
\ee

{}From Theorem \ref{t21} we conclude that the operators
\begin{equation}\label{e27}
U_{\h_2\h_1}^{\rm{\scriptscriptstyle H}}=C_{\h_2}\circ C_{\h_1}^{-1}:\,\Hh_{\h_1}
\rightarrow \Hh_{\h_2}
\end{equation}
define unitary transformations between the fibers which satisfy
\begin{equation*}
U_{\h_3\h_2}^{\rm {\scriptscriptstyle H}}\circ
U_{\h_2\h_1}^{\rm{\scriptscriptstyle H}}=U_{\h_3\h_1}^{\rm{\scriptscriptstyle H}}.
\end{equation*}
These operators correspond to the parallel
transport with respect to an Hermitian connection $\dh$ on $\Hh$.
Let $\hat{K}$ denote the set of (equivalence class of)
irreducible unitary representations of $K$. It follows from
Theorem \ref{t21} that by choosing an orthonormal basis in the vector spaces
$V_R$ of all the representations $R\in\hat{K}$ the matrix entries
$R_{ij}(\cdot)$ analytically continued to $K_\C$ form an orthogonal
basis of $\Hh_{\h}$
\begin{equation}\label{e29}
\nonumber
\left \{ R_{ij}(\cdot) \right \}_{i,j=1,\ldots , \ d_R}^{R\in\hat{K}}\subset \Hh_{\h}
\end{equation}
for all $\h >0$, where $d_R$ is the dimension of $R$, and
\begin{equation}\label{e211}
\nonumber
\int_{K_\C} \overline{R_{ij}(g)}\, R_{i'j'}(g)\, d\nu_\h(g)
= \frac {e^{\h c_R}}{d_R}\ \delta_{ii'}\delta_{jj'}
\end{equation}
where $c_R$ is the eigenvalue of $-\Delta$ on the representation
$R$. Therefore the sections of $\Hh$
\begin{equation*}
R_{ij}(\cdot)
\end{equation*}
form a global orthogonal frame and obviously the sections
\begin{equation}\label{e213}
{ F}^{R_{ij}}_\h=e^{-\h c_R/2} \, \sqrt{d_R} \, R_{ij}(\cdot)
\end{equation}
form a global orthonormal frame, that is, their norms do not
depend on $\h$. Moreover
\begin{equation}\label{e214}
U_{\h_2\h_1}^{\rm{\scriptscriptstyle H}}(R_{ij})=
e^{-(\h_2-\h_1)c_R/2}\, R_{ij} \ \Leftrightarrow  \ U_{\h_2\h_1}^{\rm{\scriptscriptstyle H}}
({ F}^{R_{ij}}_{\h_1}) = { F}^{R_{ij}}_{\h_2} \ .
\end{equation}

\begin{definition} \label{d21}
The \emph{CST connection} $\dh$ on $\Hh$ is the (unique) connection
for which the sections (\ref{e213}) are horizontal
\begin{equation}\label{e215}
\dh ({ F}^{R_{ij}} )=0.
\end{equation}
\end{definition}
The CST connection is automatically Hermitian.
{}From (\ref{e213}), $\dh$ can be equivalently defined through
\begin{equation*}
\dh_{\frac{\partial}{\partial\h}}(R_{ij})=\frac{c_R}{2}\,
R_{ij}=-\frac{\Delta_\C}{4}\, R_{ij}.
\end{equation*}

\begin{proposition}\label{p21}
The unitary parallel transport $U$ corresponding to the connection (\ref{e215}) 
is given by the unitary operators (\ref{e27}), $U=U^{\rm H}$.
\end{proposition}

\begin{proof}
Since the sections ${ F}^{R_{ij}}$ satisfy (\ref{e215})
the parallel transport for them is
\begin{equation}\label{e219}
U_{\h_2\h_1}{ F}^{R_{ij}}_{\h_1}={ F}^{R_{ij}}_{\h_2}\qquad
\forall \h_1,\h_2\in\R_+.
\end{equation}
This map of orthonormal basis extends to a unique unitary isomorphism
\begin{equation*}
U_{\h_2\h_1}:\Hh_{\h_1}\rightarrow \Hh_{\h_2}
\end{equation*}
{}From (\ref{e214}) and (\ref{e219}) we see that $U_{\h_2\h_1}$ and
$U_{\h_2\h_1}^{\rm{\scriptscriptstyle H}}$ coincide on the basis vectors and therefore
they coincide as unitary operators.
\end{proof}


\subsection{Equivalence between $\dh$ and $\dq$}
\la{s4.2}

Our main result in this section is

\begin{theorem}
\label{t41}
The quantum connection $\dq$ and the CST connection $\dh$
are equivalent in the sense that there exists  a natural unitary
isomorphism \newline $S:\Hh \to \Hq$ such that 
\be
\label{comuta}
\dq_{\frac{\partial}{\partial \h}}\circ S = S \circ \dh_{\frac{\partial}{\partial \h}}.
\ee
\end{theorem}

\begin{proof}
We start by constructing a natural unitary isomorphism between the CST bundle $\Hh$
and the quantum bundle $\Hq$. Let $\sigma, \zeta$ be two sections
of $\Hq$ of the form (\ref{ts}), (\ref{ts2}),
\ba
\nonumber
\sigma_s &=& \hat F(xe^{isY})\ e^{-\frac{s|Y|^2}{2\h_0}} \sqrt{\Omega_s}\\
\zeta_s &=& \hat G(xe^{isY})\ e^{-\frac{s|Y|^2}{2\h_0}} \sqrt{\Omega_s} .
\nonumber
\ea
{} From proposition \ref{l31} we obtain
\be
\nonumber
\langle \sigma_s, \zeta_s \rangle^Q = a_s \int_{K_\C}
\overline{ \hat F} \hat G d\nu_\h  \ ,
\ee
where $\h = s\h_0$ and $a_s$ was defined in (\ref{sal2a}).
Therefore, the bundle morphism
\be
\nonumber
\begin{array}{ccccl}
S_\h & : & \Hh_\h & \rightarrow & \Hq_s \\
     &   &  \hat F & \mapsto & \psi_s^* (\hat F)\,
 e^{-\frac{s|Y|^2}{2\h_0}}\sqrt{\frac {\Omega_s}{a_s}} ,
\end{array}
\ee
is a unitary isomorphism.
 To show (\ref{comuta}) it is sufficient
to see that the frame of horizontal sections 
$$
\left\{ F^{R_{ij}}
= e^{-\h c_R/2}R_{ij}\right\}
$$
is mapped to an horizontal frame. This follows directly from theorem \ref{th3.1}.

\end{proof}

\vskip 2cm


{\large{\bf Acknowledgements:}} We would also like to thank the referee
for the interest, the careful reading of the manuscript
and useful suggestions. The authors were partially
supported by the Center for Mathematics and its Applications, 
IST (CF, PM and JM), the Multidisciplinary Center of Astrophysics, IST (PM), 
the Center for Analysis, Geometry and Dynamical Systems, IST (JPN), and
by the Funda\c c\~ao para a Ci\^encia e a Tecnologia (FCT) through
the programs PRAXIS XXI, POCTI, FEDER, the project POCTI/33943/MAT/2000, 
and the project CERN/FIS/43717/2001.

\vspace{1cm}

\noindent \small{$\dagger$Department of Mathematics, Instituto 
Superior T\'ecnico, Av. Rovisco Pais,
1049-001 Lisboa, Portugal \\ Email: cfloren, jmourao, jpnunes@math.ist.utl.pt}

\vskip 0.2cm

\noindent \small{$\ddagger$Center for Mathematics and its Applications, Instituto Superior
T\'ecnico, Av. Rovisco Pais, 1049-001 Lisboa, Portugal \\ Email:
pmatias@math.ist.utl.pt}


\begin{thebibliography}{MMM}


\small

\baselineskip 4.5mm

\bibitem[AdPW]{AdPW} S.Axelrod, S.Della Pietra, E.Witten, ``Geometric
quantization of
Chern-Simons gauge theory'', J. Diff. Geom. {\bf 33} (1991), 787-902.




\bibitem[Fa]{Fa} G.Faltings, ``Stable $G$-bundles and projective connections'',
J. Algebraic Geom. {\bf 2} (1993) 507-568.

\bibitem[FMMN]{FMMN} C.Florentino, P.Matias, J.Mour\~ao, J.P.Nunes, ``On the BKS pairing
for K\"ahler quantizations of the cotangent bundle of a Lie group'', in preparation. 

\bibitem[FMN]{FMN} C.Florentino, J.Mour\~ao, J.P.Nunes, ``Coherent 
state transform and abelian varieties'', J. Funct. Anal. {\bf 192} (2002) 410-424;\\
``Coherent State Transforms and Vector Bundles
on Elliptic Curves'', J. Funct. Anal. {\bf 204} (2003) 355-398.

\bibitem[Ha1]{Ha1} B.C.Hall, ``The Segal-Bargmann coherent state
transform for compact Lie groups'', J.
Funct. Anal. {\bf 122} (1994), 103-151.


\bibitem[Ha2]{Ha2} B.C.Hall, ``Phase space bounds for quantum mechanics on a 
 compact Lie group'',
 Comm. Math. Phys. {\bf 184} (1997) 233-250;\\
``Harmonic analysis with respect to the heat kernel measure'',
Bull. Amer. Math. Soc. {\bf 38} (2001) 43-78.



\bibitem[Ha3]{Ha3} B.C.Hall, ``Geometric quantization and the generalized 
Segal-Bargmann transform for
Lie groups of compact type'', Comm. Math. Physics {\bf 226} (2002) 233-268;




\bibitem[Hi]{Hi} N.Hitchin, ``Flat connections and geometric quantization'', 
Comm, Math. Phys. {\bf 131}
(1990) 347-380.


\bibitem[Las]{Las} Y. Laszlo,
``Hitchin's and WZW connections are the same'',
J. Differential Geom. {\bf 49} (1998)547-576.


\bibitem[Ra]{Ra} T.R. Ramadas,
``Faltings' construction of the K-Z connection'',
Comm. Math. Phys. {\bf 196} (1998) 133-143.

\bibitem[Th1]{Th1} T.Thiemann, ``Reality Conditions Inducing Transforms 
for Quantum Gauge Field Theory
and Quantum Gravity'', Class. Quantum Grav. {\bf 13} (1996)
1383-1403.

\bibitem[Th2]{Th2} T.Thiemann, ``Gauge Field Theory Coherent States I. 
General Properties'', Class. quantum Grav.
{\bf 18}, 2025-2064.


\bibitem[Wo]{Wo} N. Woodhouse, ``Geometric Quantization'', Clarendon Press, Oxford, 1992.

\end{thebibliography}
\end{document}